\documentclass[12pt]{article}
\usepackage{amsmath,amsthm}
\usepackage{amssymb,bbm}
\usepackage[T1]{fontenc}

\newtheorem{theorem}{Theorem}[section]
\newtheorem{lemma}[theorem]{Lemma}

\newcommand{\Set}[2]{
\left\{#1\,\middle|\,#2\right\}}
\newcommand{\id}{\mathrm{id}}
\newcommand{\N}{\mathbb{N}}
\newcommand{\R}{\mathbb{R}}
\newcommand{\M}{\mathcal{M}}
\newcommand{\Pow}{\mathfrak{P}}
\newcommand{\ELL}{\overline{\ell^\infty_+}}
\newcommand{\dGH}{d_{\mathrm{GH}}}

\newcommand{\dis}{\mathop{\mathrm{dis}}\nolimits}
\newcommand{\diam}{\mathop{\mathrm{diam}}\nolimits}
\newcommand{\GHC}{\mathop{\mathrm{GHC}}\nolimits}

\newcommand{\sym}{\mathrel{\triangle}}

\begin{document}

\title{An isometric embedding of the space of non-increasing sequences with non-negative terms into the Gromov-Hausdorff space}

\author{Takuma Byakuno}
\date{}
\maketitle

\renewcommand{\thefootnote}{}
\thefootnote{2020 \emph{Mathematics Subject Classification}: Primary 30L05; Secondary 54E35, 53C23.}
\thefootnote{\emph{Key words and phrases}: Gromov-Hausdorff space, Isometric embeddings.}
\renewcommand{\thefootnote}{\arabic{footnote}}
\setcounter{footnote}{0}

\begin{abstract}
The Gromov-Hausdorff space is a certain metric space induced by the Gromov-Hausdorff convergence, and we study an isometric embedding into this space.
While the subspace of non-negative real numbers in the Euclidean space admits an isometric embedding, we prove that the subspace of $\ell^\infty$ consisting of non-increasing sequences with non-negative terms can be isometrically embedded into the Gromov-Hausdorff space.
\end{abstract}

\section{Introduction}
In this paper, we write $(\M,\dGH)$ as the {\it Gromov-Hausdorff space} for a certain set $\M$ and a metric $\dGH$.
The set $\M$ is the quotient of the class of all non-empty compact metric spaces by isometry, and the function $\dGH$, called the {\it Gromov-Hausdorff distance}, measures the difference in shape between compact metric spaces.

Convergence with respect to the Gromov–Hausdorff distance $\dGH$, called {\it Gromov–Hausdorff convergence}, plays a role in some modern geometric theories.
The Gromov-Hausdorff space is expected to possess interesting topological properties, as it is a large space endowed with a geometric convergence.
Recently, Ishiki \cite{Ishiki} proved that the Gromov-Hausdorff space is homeomorphic to the Hilbert space.

In this paper, we study isometric embeddings into the Gromov-Hausdorff space.
Iliadis, Ivanov, and Tuzhilin \cite{Tuzhilin1} showed that all finite metric spaces can be isometrically embedded into the Gromov-Hausdorff space.
In the same paper, the authors raised the question of whether every compact metric space admits such an isometric embedding.
In \cite{Byakuno}, for a fixed sequence $(r_n)_{n=0}^\infty$ with $r_n>0$ and $\sum r_n<\infty$,  we embed $\prod[0,r_n]\subset\ell^\infty$ into the Gromov-Hausdorff space and gave a partial answer.
Fukano \cite{Fukano} gave an affirmative answer to this problem.
While the case of bounded separable metric spaces arises as the next natu
ral question, in this paper we embed a particular unbounded separable metric space into the Gromov–Hausdorff space.

We write $\diam(X)$ as the diameter of a metric space $X$.
The definition of $\dGH$ yields
\[
\dGH(tX,sX)=\frac12|t-s|\diam(X)
\qquad\text{for all } X\in\M \text{ and } s,t>0, 
\]
where $tX$ is a scale transformation of $X\in\M$ by $t>0$.
Taking $X\in\M$ with $\diam(X)=2$, we obtain an isometric embedding of the half-line $[0,\infty)$ into the Gromov–Hausdorff space $\M$.
However, it was shown in \cite[Corollary 9]{Tuzhilin2} that the real line $\R$ can not admit such an isometric embedding.

As the real line $\R$ can be regarded as a subspace of $\ell^\infty$, the half-line $[0,\infty)$ has an analogue in $\ell^\infty$, namely the space of non-increasing sequences with non-negative terms.
In this paper, we prove the following.
\begin{theorem}
\label{mthm}
When the set of non-negative integers is denoted by $\N$, 
The subspace
\[
\ELL=\Set{(a_n)_{n\in\N}}{\text{$a_n\geq a_m\geq0$ for $n<m$}}
\]
of $\ell^\infty$ can be isometrically embedded into the Gromov-Hausdorff space $\M$.
In other words, there exists a map $\Pow:\ELL\to\M$ such that
\[
\dGH(\Pow(a_n)_{n\in\N},\Pow(b_n)_{n\in\N})
=\sup_{n\in\N}|a_n-b_n|
\]
for every $(a_n)_{n\in\N}, (b_n)_{n\in\N}\in\ELL$. 
\end{theorem}

\section{Preliminaries}
Let $X$ and $Y$ be sets, and $\Gamma\subset X\times Y$.
We denote by $s:\Gamma\to X$ and $t:\Gamma\to Y$ the projections onto $X$ and $Y$, respectively.
For each $e\in \Gamma$, we then have $e=(s(e),t(e))$.

Let $n\in\N$. 
we set $[n]=\{0, 1, \ldots, n-1\}$ and denote by $\Pow_n$ the power set of $[n]$.
The power set $\Pow_n$ becomes a group under the symmetric difference
\[
x\sym y = (x\cup y)\setminus(x\cap y)
=(x\setminus y) \cup (x\setminus y)
\qquad\text{for all } x,y\in\Pow_n.
\]
We equip the Cartesian product $\Pow_n\times\Pow_n$ with the product operation
\[
e_0\sym e_1=(s(e_0)\sym s(e_1), t(e_0)\sym t(e_1))
\qquad\text{for all }  e_0,e_1\in\Pow_n\times\Pow_n.
\]

\par
In the proof of Theorem~\ref{mthm}, we will equip $\Pow_n$ with a metric and calculate the Gromov-Hausdorff distance between them.

A subset $\Gamma\subset X\times Y$ is sometimes called a {\it correspondence} between $X$ and $Y$.
We define $\GHC(X, Y)$ as the set of all correspondences $\Gamma\subset X\times Y$ satisfying $s[\Gamma]=X$ and $t[\Gamma]=Y$.
Let $d_X$ and $d_Y$ denote the metric on $X$ and $Y$, respectively.
For any $\Gamma\in\GHC(X,Y)$, we define its {\it distortion} by
\[
\dis(\Gamma,d_X,d_Y)
=\sup_{e_0,e_1\in\Gamma}|d_X(s(e_0),s(e_1))-d_Y(t(e_0),t(e_1))|.
\]
The Gromov-Hausdorff distance between metric spaces $(X,d_X)$ and $(Y,d_Y)$ is then given by
\[
\dGH(X,Y)
=\frac12\inf_{\Gamma\in\GHC(X,Y)}\dis(\Gamma,d_X,d_Y).
\]
For $t>0$, we denote by $tX$ the metric space $X$ equipped with the metric multiplied by $t$.
It is straightforward to verify from the definition of $\dGH$ that
\[
\dGH(tX,tY)=t\dGH(X,Y).
\]
It is well-known that $\dGH$ defines a metric on $\M$, the quotient of the class of all non-empty compact metric spaces by isometry.
We call the metric space $(\M,\dGH)$ the {\it Gromov-Hausdorff space}.
It is known that this space is complete.

We set $\Pow_n^*=\Pow_n\setminus\{\emptyset\}$.
In the proof of Theorem~\ref{mthm}, the computation of the Gromov–Hausdorff distance produces a specific correspondence between $\Pow_n^*$ and itself.
We will then apply Hall's theorem to this correspondence.

Let $\Gamma\subset X\times Y$ be a correspondence between sets $X$ and $Y$.
For any subset $A\subset X$, we define the {\it image} of $A$ under $\Gamma$ by
\[
\Gamma[A]
=\Set{y\in Y}{(a,y)\in\Gamma,\, a\in A}
=\Set{t(e)\in Y}{e\in\Gamma,\, s(e)\in A}.
\] 
We view a map $f: X\to Y$ as the graph of $f$, which is a special correspondence
\[
\Set{(x,f(x))\in X\times Y}{x\in X}.
\]
With these conventions, Hall's marriage theorem can be stated as follows.
\begin{theorem}[Hall's marriage theorem \cite{Hall}]
Let $\Phi\subset X\times Y$ be a correspondence between finite sets $X$ and $Y$.
The following are equivalent.
\begin{enumerate}
\item
$\Phi$ includes a bijective map $\phi: X\to Y$; that is,
\[
\Set{(x,\phi(x))\in X\times Y}{x\in X}\subset\Phi.
\]
\item
Every subset $A\subset X$ satisfies the Hall condition $|A|\leq|\Phi[A]|$, where $|A|$ is the cardinarity of $A$.
\end{enumerate}
\end{theorem}

To prove Theorem~\ref{mthm}, we divide the proof into two sections: the first treats the finite-dimensional case, and the second extends the argument to the whole space.

\section{Embedding of finite-dimensional subspace}

We recall the definition of the space
\[
\ELL=\Set{(a_n)_{n\in\N}}{\text{$a_n\geq a_m\geq0$ for $n<m$}}.
\]
For $(a_n)_{n\in\N}\in\ELL$ and $x\subset\N$, we set
\[
\langle(a_n)_{n\in\N}, x\rangle=\sup_{n\in x}a_n,
\]
with the convention that
\[
\langle(a_n)_{n\in\N}, \emptyset\rangle=\sup\emptyset=0.
\]
If $x\neq\emptyset$, then $\langle(a_n)_{n\in\N}, x\rangle=a_{\min x}$.

Fix $n\in\N$.
We identify any sequence $(a_k)_{k\in[n]}$ of length $n$ with a sequence $(a_k)_{k\in\N}$ by setting $a_k=0$ for all $k\geq n$.
In this section, we embed the subspace 
\[
\R^n_+=\Set{(a_k)_{k\in[n]}}{a_0\geq a_1\geq\cdots\geq a_{n-1}>0}
\subset\ELL
\]
into the Gromov-Hausdorff space $\M$.

Let $(a_k)_{k\in[n]}\in\R^n_+$.
We define a metric on $\Pow_n$ by
\[
d(x,y)=\langle(a_k)_{k\in[n]}, x\sym y\rangle
=\max_{k\in x\sym y} a_k
\]
for $x,y\in\Pow_n$ ($x,y\subset[n]\subset\N$).
The metric space $(\Pow_n,d)$ is isometric to the $\ell^\infty$-product $\prod_{k\in[n]} [2]_k$, where each $[2]_k$ is a two-point metric space of diameter $a_k$.
Therefore, the space is compact.
We define $P_n(a_k)_{k\in[n]}\in\M$ to be this compact metric space $(\Pow_n,d)$.

To prove that $\R^n_+$ admits an isometric embedding into $\M$, it is sufficient to show that
\[
2\dGH(P_n(a_k)_{k\in[n]},P_n(b_k)_{k\in[n]})
=\max_{k\in[n]}|a_k-b_k|
\]
for all $(a_k)_{k\in[n]},(b_k)_{k\in[n]}\in\R^n_+$.

Let $d_a$ and $d_b$ denote the metrics on $P_n(a_k)_{k\in[n]}$ and $P_n(b_k)_{k\in[n]}$, respectively.
Fix $\Gamma\in\GHC(\Pow_n,\Pow_n)$.
We define a correspondence $\Phi\subset\Pow_n^*\times\Pow_n^*$ by
\[
\Phi=\Set{e_0\sym e_1}{e_0,e_1\in\Gamma,\,s(e_0)\neq s(e_1),\, t(e_0)\neq t(e_1)}.
\]
Then we obtain
\begin{align*}
\dis(\Gamma,d_a,d_b)
=&\max_{e_0,e_1\in\Gamma}|d_a(s(e_0),s(e_1))-d_b(t(e_0),t(e_1))|\\
=&\max_{e_0,e_1\in\Gamma}\Bigl|\langle(a_k)_{k\in[n]}, s(e_0)\sym s(e_1)\rangle - \langle(b_k)_{k\in[n]}, t(e_0)\sym t(e_1)\rangle\Bigr|\\
\geq&\max_{e\in\Phi}|\langle(a_k)_{k\in[n]}, s(e)\rangle - \langle(b_k)_{k\in[n]}, t(e)\rangle|.
\end{align*}

\par
By choosing $f(x)\in\Gamma[\{x\}]$ for all $x\in\Pow_n$, we obtain a map $f:\Pow_n\to\Pow_n$.
Since the inverse set $\Gamma^{-1} = \Set{(y,x)}{(x,y)\in\Gamma}$ is also a correspondence, 
we similarly obtain a map $g:\Pow_n\to\Pow_n$ satisfying $g(y)\in\Gamma^{-1}[\{y\}]$ for all $y\in\Pow_n$.
Moreover, since $f(\emptyset)\in\Gamma[\{\emptyset\}]$, we may assume that $g\circ f(\emptyset)=\emptyset$.

We can apply the following lemma to the maps $f$ and $g$ defined above.
\begin{lemma}
\label{lemma:bijection}
Let $f,g:\Pow_n\to\Pow_n$ satisfy $g\circ f(\emptyset)=\emptyset$.
Define a correspondence $\Phi^f_g\subset\Pow_n^*\times\Pow_n^*$ by
\[
\Phi^f_g=\Set{(x\sym g(y),\,f(x)\sym y)}{x,y\in\Pow_n,\, x\neq g(y),\, f(x)\neq y}.
\]
Then $\Phi^f_g$ induces a bijection $\phi:\Pow_n^*\to\Pow_n^*$; that is, $\phi\subset\Phi^f_g$.
\end{lemma}
\begin{proof}
By Hall's marriage theorem, it is sufficient to show that $|A|\leq|\Phi^f_g[A]|$ for $A\subset\Pow_n^*$.
Define a map $\iota: A\times\Pow_n\to\Pow_n\times(\Phi^f_g[A]\cup\{\emptyset\})$ by setting
\[
\iota(a,y)=(a\sym g(y),\,f(a\sym g(y))\sym y).
\]
Since $A\subset \Pow_n^*=\Pow_n\setminus\{\emptyset\}$, we have 
\[
\iota(a,y)\neq(\emptyset,\emptyset)
\qquad\text{for all } (a,y)\in A\times\Pow_n.
\]
Moreover, by using $g\circ f(\emptyset)=\emptyset$, it is easy to show that $\iota$ is injective.
Therefore, from $|\Pow_n|=2^n$, $|\Phi^f_g[A]\cup\{\emptyset\}|=|\Phi^f_g[A]|+1$, we obtain
\[
2^n|A|=|A\times\Pow_n|
\leq|(\Pow_n\times(\Phi^f_g[A]\cup\{\emptyset\}))\setminus\{(\emptyset,\emptyset)\}|
=2^n(|\Phi^f_g[A]|+1)-1.
\]
Rearranging this inequality yields $|A|\leq|\Phi^f_g[A]|$.
\end{proof}

\par
Since $\Phi^f_g\subset\Phi$,  $\Phi$ includes a bijection $\phi:\Pow_n^*\to\Pow_n^*$.
Thus, we obtain
\begin{align*}
\dis(\Gamma,d_a,d_b)
\geq&\sup_{e\in\Phi}|\langle(a_k)_{k\in[n]}, s(e)\rangle - \langle(b_k)_{k\in[n]}, t(e)\rangle|\\
\geq&\sup_{x\in\Pow_n^*}|\langle(a_k)_{k\in[n]}, x\rangle - \langle(b_k)_{k\in[n]}, \phi(x)\rangle|
\geq&\sup_{x\in\Pow_n^*}|a_{\min x} - b_{\min\phi(x)}|.
\end{align*}

It is straightforward to prove the following lemma.
\begin{lemma}
\label{lemma:undo}
Let $c_0,c_1,d_0,d_1\geq0$ satisfy $c_0\geq c_1$ and $d_0\geq d_1$.
Then we have
\[
\max\{|c_0-d_0|,|c_1-d_1|\}\leq\max\{|c_0-d_1|,|c_1-d_0|\}.
\]
\end{lemma}
\begin{lemma}
\label{lemma:sort}
Let $M\in\N$ and $\sigma:[M]\to[M]$ a bijection map.
For $(c_k)_{k\in[M]},(d_k)_{k\in[M]}\in\R^M_+$, the  inequality
\[
\sup_{k\in[M]}|c_k-d_k|
\leq\sup_{k\in[M]}|c_k-d_{\sigma(k)}|
\]
holds.
\end{lemma}
\begin{proof}
We set $\sigma_0=\sigma$.
For the induction step, assume that a bijection $\sigma_m:[M]\to [M]$ ($m\in[M]$) satisfies $\sigma_m|_{[m]}=\id_{[m]}$ ($[0]=\emptyset$) and
\[
\sup_{k\in[M]}|c_k-d_{\sigma_m(k)}|
\leq\sup_{k\in[M]}|c_k-d_{\sigma(k)}|.
\]
Using the transposition $(m\ \ \sigma_m(m)): [M]\to[M]$, we define a bijection $\sigma_{m+1}: [M]\to[M]$ by
\[
\sigma_{m+1}=(m\ \ \sigma_m(m))\circ\sigma_m.
\]
It is obvious to see $\sigma_{m+1}|_{[m+1]}=\id_{[m+1]}$.
Since $\sigma_m|_{[m]}=\id_{[m]}$, we have $m\leq\min\{\sigma_m(m),\sigma_m^{-1}(m)\}$.
Hence $c_m\geq c_{\sigma_m^{-1}(m)}$ and $d_m\geq d_{\sigma_m(m)}$.
Lemma~\ref{lemma:undo} yields
\[
\max\{|c_m-d_m|,|c_{\sigma_m^{-1}(m)}-d_{\sigma_m(m)}|\}\leq\max\{|c_m-d_{\sigma_m(m)}|,|c_{\sigma_m^{-1}(m)}-d_m|\}.
\]
Combining this with the induction hypothesis, we obtain
\[
\sup_{k\in[M]}|c_k-d_{\sigma_{m+1}(k)}|
\leq\sup_{k\in[M]}|c_k-d_{\sigma_m(k)}|
\leq\sup_{k\in[M]}|c_k-d_{\sigma(k)}|.
\]
Thus the recursive construction works.
By $\sigma_M=\id_{[M]}$, the desired inequality follows.
\end{proof}

\par
Since $\Pow_n^*$ is finite,  we may choose a bijection $i:\Pow_n^*\to[M]$ ($M=|\Pow_n^*|$) such that \[
\min i^{-1}(m_0)<\min i^{-1}(m_1)
\qquad
\text{whenever } m_0< m_1.
\]
Define sequences $(c_k)_{k\in[M]}\in\R^M_+$ and $(d_k)_{k\in[M]}\in\R^M_+$ by
\[
c_k=a_{\min i^{-1}(k)}
,\qquad
d_k=b_{\min i^{-1}(k)}.
\]
Applying Lemma~\ref{lemma:sort}, we obtain
\begin{align*}
\dis(\Gamma,d_a,d_b)
\geq&\sup_{x\in\Pow_n^*}|a_{\min x} - b_{\min\phi(x)}|
=\sup_{k\in[M]}|c_k-d_{i\circ\phi\circ i^{-1}(k)}|\\
\geq&\sup_{k\in[M]}|c_k-d_k|
=\sup_{x\in\Pow_n^*}|a_{\min x} - b_{\min x}|.
\end{align*}
Hence, by restricting to one-point sets $x\in\Pow_n$, we have
\[
\dis(\Gamma,d_a,d_b)\geq\sup_{k\in[n]}|a_k - b_k|.
\]

Collecting the above computations, we obtain
\[
2\dGH(P_n(a_k)_{k\in[n]},P_n(b_k)_{k\in[n]})
=\inf_{\Gamma\in\GHC(\Pow_n,\Pow_n)}\dis(\Gamma,d_a,d_b)\geq\sup_{k\in[n]}|a_k - b_k|.
\]
Conversely, since $\id_{\Pow_n}\in\GHC(\Pow_n,\Pow_n)$ and 
\[
\dis(\id_{\Pow_n},d_a,d_b)=\sup_{k\in[n]}|a_k - b_k|,
\]
we obtain the reverse inequality.
Thus the desired identity follows.

\section{Extension to the whole space}

We recall the notation from the previous section.
We proved that
\[
\R^n_+ =\Set{(a_k)_{k\in[n]}}{a_0\geq a_1\geq\cdots\geq a_{n-1}>0}
\subset\ELL
\]
admits an isometric embedding into $\M$ for $n\in\N$.
For all $(a_n)_{n\in\N}\in\ELL$ and $x\subset\N$, we define
\[
\langle(a_n)_{n\in\N}, x\rangle=\sup_{n\in x}a_n.
\]
Fix $n\in\N$.
Then the metric on $P_n(a_k)_{k\in[n]}=\Pow_n$  is given by
\[
d(x,y)=\langle (a_k)_{k\in[n]}, x\sym y \rangle
\qquad x,y\in\Pow_n.
\]
Moreover,  we have
\[
2\dGH(P_n(a_k)_{k\in[n]},P_n(b_k)_{k\in[n]})
=\max_{k\in[n]}|a_k-b_k|.
\]
Hence, the map $2P_n:\R^n_+\to\M$ is an isometric embedding.

By passing to the completion, this embedding $2P_n:\R^n_+\to\M$ extends uniquely to an isometric embedding $2P_n:\overline{\R^n_+}\to\M$, where
\[
\overline{\R^n_+}=\Set{(a_k)_{k\in[n]}}{a_0\geq a_1\geq\cdots\geq a_{n-1}\geq0}
\subset\ELL.
\]
\begin{lemma}
\label{lemma:rank}
Let $n_0,n_1\in\N$ with $n_0<n_1$, and let $(a_k)_{k\in[n_0]}\in\R^{n_0}_+\subset\ELL$.
Then $(a_k)_{k\in[n_1]}\in\overline{\R^{n_1}_+}$, and the metric space $P_{n_1}(a_k)_{k\in[n_1]}$ is isometric to $P_{n_0}(a_k)_{k\in[n_0]}$.
\end{lemma}
\begin{proof}
For $l\in\N$, we define a sequence $(a^l_k)_{k\in[n_1]}$ by setting
\[
a^l_k=\begin{cases}
a_k&\mathrm{if:}\,k\in[n_0],\\
2^{-l}&\mathrm{if:}\,k\in[n_1]\setminus[n_0].
\end{cases}
\]
Then, we have $(a^l_k)_{k\in[n_1]}\in\R^{n_1}_+$ for all sufficiently large $l$ and $(a^l_k)_{k\in[n_1]}\to(a_k)_{k\in[n_1]}$ as $l\to\infty$. 
It is straightforward to see that the correspondence 
\[
\Gamma=\Set{(y\cap[n_0],y)}{y\in\Pow_{n_1}}\subset\Pow_{n_0}\times\Pow_{n_1}
\]
belongs to $\GHC(\Pow_{n_0},\Pow_{n_1})$.

Fix $l\in\N$.
Note that $P_{n_0}(a_k)_{k\in[n_0]}$ is a subspace of $P_{n_1}(a^l_k)_{k\in[n_1]}$ by the definition of $(a^l_k)_{k\in[n_1]}$.
Denoting by $d$ the metrics on $P_{n_1}(a^l_k)_{k\in[n_1]}$, the definition of $\dGH$ yields
\[
\dGH(P_{n_0}(a_k)_{k\in[n_0]},P_{n_1}(a^l_k)_{k\in[n_1]})
\leq\frac12\dis(\Gamma,d,d)
\]
Since $y\sym(y\cap[n_0])=y\setminus[n_0]$ for $y\in\Pow_{n_1}$, 
the triangle inequality gives 
\begin{align*}
\dis(\Gamma,d,d)
=&\sup_{y_0,y_1\in\Pow_{n_1}}|d(y_0\cap[n_0],y_1\cap[n_0])- d(y_0,y_1)|\\
\leq&2\sup_{y\in\Pow_{n_1}}d(y,y\cap[n_0])
=2\sup_{y\in\Pow_{n_1}}\langle(a^l_k)_{k\in[n_1]},y\setminus[n_0]\rangle
=\frac{2}{2^l}
\end{align*}
Thus, we have $\dGH(P_{n_0}(a_k)_{k\in[n_0]},P_{n_1}(a^l_k)_{k\in[n_1]})\leq2^{-l}$.
Letting $l\to\infty$, we obtain $\dGH(P_{n_0}(a_k)_{k\in[n_0]},P_{n_1}(a_k)_{k\in[n_1]})=0$, and therefore the two space are isometric.
\end{proof}

\par
Let $(a_n)_{n\in\N}\in\ELL$. 
For $n_0,n_1\in\N$ with $n_0<n_1$, 
we have $(a_k)_{k\in[n_0]}\in\overline{\R^{n_0}_+}$ and $(a_k)_{k\in[n_1]}\in\overline{\R^{n_1}_+}$.
Lemma~\ref{lemma:rank} implies
\[
2\dGH(P_{n_0}(a_k)_{k\in[n_0]},P_{n_1}(a_k)_{k\in[n_1]})
=\max_{k\in[n_1]\setminus[n_0]}|0-a_k|
=a_{n_0}.
\]
Since $(a_n)_{n\in\N}$ is bounded and non-increasing, it converges.
Therefore, the sequence of compact metric spaces $(P_n(a_k)_{k\in[n]})_{n\in\N}$ is a Cauchy  sequence in $\M$ and we obtain its limit $P(a_n)_{n\in\N}$.

We now show that $2P:\ELL\to\M$ is an isometric embedding.
Let $(a_n)_{n\in\N},(b_n)_{n\in\N}\in\ELL$. 
For each $n\in\N$, since $(a_k)_{k\in[n]},(b_k)_{k\in[n]}\in\overline{\R^n}_+$, we obtain
\[
2\dGH(P_n(a_k)_{k\in[n]},P_n(b_k)_{k\in[n]})=\max_{k\in[n]}|a_k-b_k|.
\]
Letting $n\to\infty$, we have hence
\[
2\dGH(P(a_n)_{n\in\N},P(b_n)_{n\in\N})=\sup_{n\in\N}|a_n-b_n|.
\]
Consequently, $\ELL$ can be isometrically embedded into $\M$.

\section{Acknowledgments}
This research was financially supported by the Kansai University Grant-in-Aid for research progress in graduate course, 2026.
The author acknowledges this support.
He is grateful to Dr.~Yoshito~Ishiki and Professor~Atsushi~Fujioka for their encouragement.
He also thanks Professor~Kohji~Yanagawa and GPT-6 for their helpful comments on the proof of Lemma~\ref{lemma:bijection}.
The author sincerely appreciates their support.

Takuma Byakuno: the Graduate School of Science and Engineering, Kansai University, 3-3-35 Suita, 564-8680 Osaka, Japan, E-mail: k241344@kansai-u.ac.jp; tatatagobyakuno0822@gmail.com


\begin{thebibliography}{9}
\bibitem{Tuzhilin2}
S.~Borzov, A.~Ivanov, A.Tuzhilin, 
\newblock{\em Geometry of the Gromov-Hausdorff distance on the class of all metric spaces},
\newblock Mathematics, 2022, Volume 213, Issue 5, 641–658.

\bibitem{Byakuno}
T.~Byakuno, 
\newblock{\em A sequence of compact metric spaces and an isometric embedding into the Gromov-Hausdorff space},
\newblock Topology and its Applications, {\bf 384}, 2026, DOI: 10.1016/j.topol.2026.109790.

\bibitem{Fukano}
R.~Fukano,
\newblock {\em Every Compact Metric Space Is Isometrically Embeddable into the Gromov-Hausdorff Space},
\newblock 2026. arXiv:2609.21469.

\bibitem{Hall}
P.~Hall,
\newblock{\em On representatives of subsets},
\newblock Journal of the London Mathematical Society, {\bf 10}, 1935.

\bibitem{Tuzhilin1}
S.~Iliadis, A.~Ivanov and A.~Tuzhilin.
\newblock {\em Local structure of Gromov-Hausdorff space, and isometric embeddings of finite metric spaces into this space}, 
\newblock Topology and its Applications, {\bf 221}, 2017, pp.393--398. 
DOI: 10.1016/j.topol.2017.02.050.

\bibitem{Ishiki}
Y.~Ishiki, 
\newblock {\em The topology of Gromov--Hausdorff space}, 
\newblock 2026. arXiv: 2609.09639.
\end{thebibliography}
\end{document}